\documentclass{amsart}
\usepackage{amssymb}
\usepackage{amsfonts}
\usepackage{amsmath}
\usepackage{graphicx}
\usepackage[pagebackref]{hyperref}
\theoremstyle{plain}
\newtheorem{thm}{Theorem}[section]
\newtheorem{corollary}[thm]{Corollary}
\newtheorem{lemma}[thm]{Lemma}
\newtheorem{proposition}[thm]{Proposition}

\theoremstyle{definition}

\theoremstyle{remark}

\newtheorem{example}[thm]{Example}

\newtheorem{claim}{\rm Claim}
\numberwithin{equation}{section}
\newcommand{\field}[1]{\mathbb{#1}}

\newcommand{\Q}{\field{Q}}

\def\1{{\rm (1)}}
\def\2{{\rm (2)}}
\def\3{{\rm (3)}}
\def\4{{\rm (4)}}
\def\5{{\rm (5)}}
\def\i{{\rm (i) }}
\def\ii{{\rm (ii) }}
\def\iii{{\rm (iii) }}
\def\iv{{\rm (iv) }}

\begin{document}

%%%%%%%%%%%%%%%%%%%%%%%%%%%%%%%%%%%%%%%%%%%%%%%%%%%%%%%%%
%%%%%%%%%%%%%%%%%%%%%%%%%%%%%%%%%%%%%%%%%%%%%%%%%%%%%%%%%
\title[Duals and Transforms of Ideals in PVMDs]
{Duals and Transforms of Ideals in PVMDs}

%%%%%%%%%%%%%%%%%%%%%%%%%%%%%%%%%%%%%%%%%%%%%%%%%%%%%%%%%
%%%%%%%%%%%%%%%%%%%%%%%%%%%%%%%%%%%%%%%%%%%%%%%%%%%%%%%%%
\author{A. BenObaid}

\address{Department of Mathematical Sciences,
King Fahd University of Petroleum \& Minerals,
Dhahran 31261, Saudi Arabia}

\email{g200502650@kfupm.edu.sa}

\author{A. Mimouni}

\address{Department of Mathematical Sciences,
King Fahd University of Petroleum \& Minerals, P. O. Box 278,
Dhahran 31261, Saudi Arabia}

\email{amimouni@kfupm.edu.sa}

\thanks{This work is supported by KFUPM}


%%%%%%%%%%%%%%%%%%%%%%%%%%%%%%%%%%%%%%%%%%%%%%%%%%%%%%%%%%%%%%%%%%%%%%%%%%%%%%%

%%%%%%%%%%%%%%%%%%%%%%%%%%%%%%%%%%%%%%%%%%%%%%%%%%%%%%%%%%%%%%%%%
\begin{abstract}
In this paper we study when the dual of a $t$-ideal in a $PVMD$ is a
ring? and we treat the question when it coincides with its endomorphism ring. We also study particular classes of overrings
of PVMDs. Specially, we investigate the Nagata transform
and the endomorphism ring of ideals in PVMDs in an attempt
to establish analogues for well-known results on overrings of
Pr\"ufer domains.
\end{abstract}
\maketitle
%%%%%%%%%%%%%%%%%%%%%%%%%%%%%%%%%%%%%%%%%%%%%%%%%%%%%%%%%%%%%%%%%%%%%%%%%%%%%%%
\section{Introduction}
Let $R$ be an integral domain and $K$ its quotient
field. For a nonzero fractional ideals $I$ and $J$ of $R$, we
define the fractional ideal $(I:J)=\{x\in K| xJ\subseteq I\}$. We
denote $(R:I)$ by $I^{-1}$ and we call it the dual of $I$ since it is isomorphic, as an $R$-module, to
$Hom_{R}(I, R)$. The Nagata transform (or ideal transform) of $I$ is defined as
$T(I)=\bigcup^{\infty}_{n=1} (R:I^n)$ and the Kaplansky transform of $I$ is
defined as $\Omega(I)=\{u\in K:ua^{n(a)}\in R,\  \hbox{a is an arbitrary element in I and}\
n(a)\ \hbox{some positive integer}\}$. The zero cohomology of $I$ over $R$ is defined by
$R^I=\bigcup^{\infty}_{n=1} (I^n:I^n)$.  It is clear that $(I:I)\subseteq
R^I\subseteq T(I)\subseteq\Omega(I)$ and  $(I:I)\subseteq
I^{-1}\subseteq T(I)\subseteq\Omega(I)$. Also we notice that $\Omega(I)$ is
a variant of the Nagata transform $T(I)$, and useful in the
case when $I$ is not finitely generated, but if $I$ is a finitely
generated ideal of $R$, then $\Omega(I)=T(I)$. It is worthwhile noting that
$\Omega(I),\ T(I),\ (I:I)$ and $R^I$ are overrings of $R$ for each
ideal $I$ in a domain $R$, while $I^{-1}$ is not, in general, a ring. Moreover, $(I:I)$ is the largest
subring of $K$ in which $I$ is an ideal and it is isomorphic to the endomorphism ring of $I$.\\

In 1968, Brewer \cite{4} proved a representation theorem for the Nagata transform $T(I)$,
when $I$ is a finitely generated ideal (which coincides in this case with $\Omega(I)$) and in 1974, Kaplansky \cite{27} gave the  complete description of the Kaplansky transform $\Omega(I)$ for
each ideal $I$ in an integral domain $R$. He proved that ``{\em if $I$ is a
nonzero ideal of $R$, then $\Omega(I)=\bigcap R_P$, where $P$ varies
over the set of prime ideals that do not contain $I$}" (this result was also obtained independently by Hays \cite{21}).
In \cite[Exercise 11, page 331]{16} Gilmer
described $T(I)$ of an ideal $I$ which is contained
in a finite number of minimal prime ideals in a Pr\"ufer domain $R$, specifically, ``{\em let
$R$ be a Pr\"ufer domain, $I$ a nonzero ideal of $R$,
$\{P_{\alpha}\}$ the set of minimal prime ideals of $I$, and
$\{M_{\beta}\}$ the set of maximal ideals that do not contain $I$.
Then $T(I)\subseteq(\bigcap R_{Q_{\alpha}})\cap(\bigcap
R_{M_{\beta}})$, where $Q_{\alpha}$ is the unique prime ideal
determined by $\bigcap^{\infty}_{n=1}
I^nR_{P_{\alpha}}=Q_{\alpha}R_{P_{\alpha}}$. Moreover, if the set
$\{P_{\alpha}\}$ is finite, equality holds}" (see also \cite[Theorem 3.2.5]{13}). In \cite{13}, Fontana, Huckaba and Papick described some relations
between the above overrings in the case of Pr\"ufer domains. For instance, they
showed that ``{\em if $P$ is a nonzero non-invertible prime ideal
of a Pr\"ufer domain $R$, then there is no proper overring between $P^{-1}$ and
$\Omega(P)$}" (\cite[Theorem 3.3.7]{13}). In 1986, Houston \cite{19} studied the divisorial prime ideals in
PVMDs, and among others, he proved that ``{\em if $P$ is a nonzero,
non-$t$-maximal $t$-prime ideal of a PVMD $R$, then $P^{-1}=R_P\cap
\mathcal{C}_t(I)$, where $\mathcal{C}_t(I)=\displaystyle\bigcap_{I\nsubseteqq
M_{\beta}\in Max_t(R)} R_{M_{\beta}}$, and $T(P)=R_{P_0}\cap
\mathcal{C}_t(I)$, where $P_0=(\bigcap_nP^nR_P)\cap R$}" (\cite[Proposition 1.1 and Proposition 1.5]{19}).\\

Many papers in the literature deal with the fractional ideal
$I^{-1}$. The main problem is to examine settings in which $I^{-1}$
is a ring. In 1982, Huckaba and Papick \cite{22} stated the following: ``{\em
let $R$ be a Pr\"ufer domain, $I$ a nonzero ideal of $R$,
$\{P_{\alpha}\}$ the set of minimal prime ideals of $I$, and
$\{M_{\beta}\}$ the set of maximal ideals that do not contain $I$.
Then $I^{-1}\supseteq(\bigcap R_{P_{\alpha}})\cap(\bigcap
R_{M_{\beta}})$. If $I^{-1}$ is a ring, equality holds}" (\cite[Theorem 3.2 and Lemma 3.3]{22}).  They also
proved that ``{\em for a radical ideal $I$ of a Pr\"ufer domain
$R$, let $\{P_{\alpha}\}$ be the set of minimal prime ideals of $I$ and
assume that $\bigcap P_{\alpha}$ is irredundant. Then $I^{-1}$ is a
subring of $K$ if and only if for each $\alpha$, $P_{\alpha}$ is not
invertible}" (\cite[Theorem 3.8]{22}). In \cite{18}, Heinzer and Papick gave a necessary and sufficient condition
for $I^{-1}$, when it is a ring, to collapse with $(I:I)$ for an ideal $I$ in
a Pr\"ufer domain with Noetherian spectrum . Namely, they proved that ``{\em for a
Pr\"ufer domain $R$ with $Spec(R)$ Noethrian, let $I$ be a nonzero
ideal of $R$ and assume that $I^{-1}$ is a ring. Then $I^{-1}=(I:I)$ if and
only if $I=\sqrt{I}$ (i.e. $I$ is a radical ideal) if and only if the minimal prime ideals of $I$
in $(I:I)$ are all maximal ideals}" (\cite[Theorem 2.5]{18}). In 1993, Fontana, Huckaba, Papick and Roitman \cite{14} studied the
endomorphism ring of an ideal in a Pr\"ufer domain. One of their main results asserted that ``{\em for a nonzero ideal $I$ of a Pr\"ufer domain $R$, let $\{Q_{\alpha}\}$ be the set of maximal prime ideals of $\mathcal{Z}(R,I)$ and $\{M_{\beta}\}$ be the set
of maximal ideals that do not contain $I$. Then
$(I:I)\supseteq(\bigcap R_{Q_{\alpha}})\cap(\bigcap R_{M_{\beta}})$.
Moreover, if $R$ is a $QR$-domain, equality holds}" (\cite[Theorem 4.11 and Corollary 4.4]{14}). Finally in 2000, Houston, Kabbaj, Lucas and Mimouni \cite{20}, gave several characterizations of when $I^{-1}$ is a ring for a nonzero ideal $I$ in an integrally
closed domain. For instance they generalized \cite[Theorem 4.11]{14} to the $PVMD$'s case. Namely they proved that ``{\em if $I$ is an ideal of a $PVMD$ with no embedded primes, then $I^{-1}$ is a ring if and only if
$I^{-1}=(I:I)=R_{\mathfrak{N}}\cap \mathcal{C}_t(I)$, where
$\mathfrak{N}$ the complement in $R$ of the set of zero divisors on
$R/I$}" (\cite[Theorem 4.7]{20}).\\

The purpose of this paper is to continue the investigation of when the
dual of an ideal in a $PVMD$ is a ring and when it coincides with its endomorphism ring.
We also aim at giving a full description of the Nagata and Kaplansky transforms of ideals
in a $PVMD$, seeking generalizations or $t$-analogues of well-known results.\\

In Section~\ref{sec:2}, we deal with the dual of a $t$-ideal in a $PVMD$. We give a
generalization of the above mentioned results of Huckaba-Papick and Heinzer-Papick.
Precisely, we prove that ``{\em for a radical $t$-ideal $I$ of a $PVMD$ $R$, let $\{P_{\alpha}\}$ be the
set of minimal prime ideals of $I$ and assume that $\bigcap
P_{\alpha}$ is irredundant. Then $I^{-1}$ is a subring of $K$ if and
only if $P_{\alpha}$ is not $t$-invertible for each $\alpha$}" (Theorem~\ref{sec:2.3}). We Also prove that
``{\em if $R$ is a $PVMD$ with $Spec_t(R)$ Noethrian, and $I$ is a
$t$-ideal of $R$ such that $I^{-1}$ is a ring, then $I^{-1}=(I:I)$ if and
only if $I=\sqrt{I}$ if and only if the minimal prime ideals of $I$ in $(I:I)$ are all $t$-maximal ideals}" (Theorem~\ref{sec:2.5}).
In the particular case where $R$ is a Pr\"ufer domain we obtain the pre-mentioned results of Huckaba-Papick and
Heinzer-Papick simply by remarking that a Pr\"ufer domain is just a $PVMD$ in which the $t$-operation is trivial, that is, $t=d$. We close this section with a  description of the endomorphism ring of a $t$-ideal in a $tQR$-domain. Particularly we give a generalization of a well-known result by Fontana et al., \cite[Corrollary 4.4 and Theorem 4.11]{14}, that is, ``{\em let $I$ be a $t$-ideal of a $PVMD$ $R$, $\{Q_{\alpha}\}$
be the set of all maximal prime ideals of $Z(R,I)$ and $\{M_{\beta}\}$ be the set of $t$-maximal ideals of $R$ that do not
contain $I$. Then $(I:I)\supseteq(\bigcap R_{Q_{\alpha}})\cap(\bigcap R_{M_{\beta}})$, and if $R$ is a $tQR$-domain then the equality holds}" (Theorem~\ref{sec:2.13}).\\

Section~\ref{sec:3} deals with Kaplansky and Nagata transforms of an ideal in a $PVMD$. Our aim is to give the $t$-analogues for many
results of Fontana-Huckaba-Papick \cite[Section 3.3]{13} for $t$-linked overrings of $PVMDs$.  Our first main theorem generalizes \cite[Theorem 3.3.7]{13} to the case of $t$-prime ideals in a $PVMD$. For instance we prove that ``{\em if $P$ is a non-$t$-invertible $t$-prime ideal of a $PVMD$ $R$, then there is no proper overring between $P^{-1}$ and $\Omega(P)$}" (Theorem~\ref{sec:3.2}). The second main theorem is a satisfactory $t$-analogue for \cite[Theorem 3.3.10]{13}, that is, ``{\em let $R$ be a $PVMD$ and $P$ a $t$-prime ideal of $R$. Then $T(P)\subsetneqq\Omega(P)$ if and only if $T(P)=R_P\cap\Omega(P)$ and $\Omega(P)\nsubseteq R_P$. Moreover, $(P\Omega(P))_{t_1}=\Omega(P)$ if and only if $\Omega(P)\nsubseteq R_P$ if and only if $P=\sqrt{I}$ for some $t$-invertible ideal}" (Theorem~\ref{sec:3.7}). Other applications of the obtained results are given.\\

Throughout this paper $R$ is an integral domain with quotient field $K$.
By a fractional ideal, we mean a nonzero $R$-submodule $I$ of of $K$ such that $dI\subseteq R$ for some nonzero element $d$ of $R$ and by a proper ideal we mean a nonzero ideal $I$ such that $I\subsetneq R$. Recall that for a fractional ideal $I$ of $R$, the $v$-closure of $I$ is the
fractional ideal $I_{v}=(I^{-1})^{-1}$ and the $t$-closure of $I$ is the ideal $I_{t}=\bigcup J_{v}$, where $J$ ranges over the set of all finitely generated subideals of $I$. A fractional ideal $I$ is said to be a $v$-ideal (or divisorial) (resp. $t$-ideal, resp. $t$-invertible) if $I=I_{v}$ (resp. $I=I_{t}$, resp. $(II^{-1})_{t}=R$), and a domain $R$ is said to be a $PVMD$ (for Pr\"ufer $v$-multiplication domain) if every nonzero finitely generated ideal is $t$-invertible (equivalently, $R_{M}$ is a valuation domain for every $t$-maximal ideal $M$ of $R$). For more details on the $v$- and $t$-operations, we refer the reader to \cite[section 32]{16}. Also it is worth to note that many of our results are inspired from the Pr\"ufer case, and some proofs are dense and use a lot of technics of the $t$-operation. We are grateful to the huge work on the $t$-move (from Pr\"ufer to $PVMD$) done during the last decades.

%%%%%%%%%%%%%%%%%%%%%%%%%%%%%%%%%%%%%%%%%%%%%%%%%%%%%%%%%%%%%%%%%%%%%%%%%%%%%%%
%%%%%%%%%%%%%%%%%%%%%%%%%%%%%%%%%%%%%%%%%%%%%%%%%%%%%%%%%%%%%%%%%%%%%%%%%%%%%%%
\section{Duals of ideals in a $PVMD$}\label{sec:2}
%%%%%%%%%%%%%%%%%%%%%%%%%%%%%%%%%%%%%%%%%%%%%%%%%%%%%%%%%%%%%
%%%%%%%%%%%%%%%% Lemma 1.1 %%%%%%%%%%%%%%%%%%%%%%%%%%%%%%%%%%
We start this section by noticing that for a fractional ideal $I$ of a
domain $R$, $I^{-1}=(I_{t})^{-1}=(I_{v})^{-1}$, $I$ is $t$-invertible if and only if $I_{t}$
is $t$-invertible and if $I_{t}=R$, then $I^{-1}=(I:I)=R$. In this regard,
we will focus on the case where $I$ is a proper $t$-ideal of $R$.\\
Before giving the first main theorem of this section, we begin with the following two results
on necessary and sufficient conditions for $I^{-1}$ to be a ring. The first one is a generalization of \cite[Lemma 2.0]{22}
(since invertible ideals are $t$-invertible $t$-ideals) and the second one is a $t$-analogue of \cite[Proposition 2.2]{20}.\\

\begin{lemma}\label{sec:2.1}
Let $R$ be a domain and $I$ a $t$-ideal of $R$. If $I$ is $t$-invertible, then $I^{-1}$ is not a ring.
\end{lemma}
%%%%%%%%%%%%%%%%%%%%%%%%%%%%%%%%%%%%%%%%%%%%%%%%%%%%%%%%%
\begin{proof}
Deny, assume that $I^{-1}$ is a ring. Let $M$ be a $t$-maximal ideal
of $R$ containing $I$. Since $I$ is $t$-invertible, then
$II^{-1}$ is not contained in any $t$-maximal ideal of $R$. Hence
$(II^{-1})_M=R_M$. So $IR_M$ is an invertible ideal of $R_M$ and
hence principal. Since $I$ is $t$-invertible, then $I$ is
$v$-finite. Hence there is a finitely generated ideal $A$ of $R$
such that $A\subseteq I$ and $I=A_t=A_v$. Since $A$ is a finitely
generated ideal of $R$, by \cite[Lemma 4]{28}, $(AR_M)_{v_1}=(A_vR_M)_{v_1}$, where $v_1$ is the $v$-operation with
respect to $R_M$. So
$(IR_M)^{-1}=(A_vR_M)^{-1}=(AR_M)^{-1}=A^{-1}R_M=(A_v)^{-1}R_M=I^{-1}R_M$.
Since $I^{-1}$ is a ring, $(IR_M)^{-1}$ is also a ring, which contradicts
the fact that $IR_M$ is principal in $R_M$.
\end{proof}
%%%%%%%%%%%%%%%%%%%%%%%%%%%%%%%%%%%%%%%%%%%%%%%%%%%%%%%%%%%%%%%%
%%%%%%%%%%%%%%%%%%%%%%%%% Corollary 1.2 %%%%%%%%%%%%%%%%%%%%%%%%
\begin{corollary}\label{sec:2.2}
Let $I$ be a $t$-ideal of a domain $R$. Then $I^{-1}$ is a ring if and only if $I$ is not $t$-invertible and $(M:I)$ is a ring for each $t$-maximal ideal $M\supseteq I$ of $R$.
\end{corollary}
%%%%%%%%%%%%%%%%%%%%%%%%%%%%%%%%%%%%%%%%%%%%%%%%%%%%%%%%%
\begin{proof}
If $I^{-1}$ is a ring, then $I$ is not $t$-invertible by Lemma~\ref{sec:2.1}. By
\cite[Proposition 2.1]{20}, $(M:I)$ is a ring for each $t$-maximal
ideal $M$ containing $I$. Conversely, if $I$ is not $t$-invertible, then
$II^{-1}\subseteq M$ for some $t$-maximal ideal $M$ of $R$ and hence
$I^{-1}=(M:I)$. So $I^{-1}$ is a ring.
\end{proof}
%%%%%%%%%%%%%%%%%%%%%%%%%%%%%%%%%%%%%%%%%%%%%%%%%%%%%%%%%%%%%%%%%
%%%%%%%%%%%%%%%%%%%%%%%%%%%%%%%%%%%%%%%%%%%%%%%%%%%%%%%%%%%%%%%%%
Now, we turn our attention to the duals of ideals in a $PVMD$. Our approach is
similar to that one of Huckaba-Papick  done in \cite{22} for Pr\"ufer domains.
Let $R$ be a $PVMD$. We divide $Spec_t(R)$, that is, the set of all nonzero $t$-prime ideals of $R$, into three disjoint sets:\\
$S_1=\{P\in Spec_t(R):P\ \hbox{is}\ t\ \hbox{-invertible}\}$\\
$S_2=\{P\in Spec_t(R):P\ \hbox{is a non-}\ t\ \hbox{-invertible}\ t\ \hbox{-maximal ideal and}\ PR_P\  \hbox{is principal}\}$\\
$S_3=\{P\in Spec_t(R):P\not\in S_1\cup S_2\}$. Our first main theorem is a generalization of \cite[Theorem 3.8]{22} to $PVMDs$.\\

%%%%%%%%%%%%%%%%%%%%%%%%%%%%%%%%%%%%%%%%%%%%%%%%%%%%%%%%%%%%%%%%%
%%%%%%%%%%%%%%%%%%%%%% ProPosition 1.3 %%%%%%%%%%%%%%%%%%%%%%%%%%
%%%%%%%%%%%%%%%%%%%%%%%%%%%%%%%%%%%%%%%%%%%%%%%%%%%%%%%%%%%%%%%%%
\begin{thm}\label{sec:2.3}
Let $I$ be a radical $t$-ideal of a $PVMD$ $R$, $\{P_{\alpha}\}$ the
set of all minimal prime ideals of $I$ and assume that $\bigcap
P_{\alpha}$ is irredundant. Then $I^{-1}$ is a subring of $K$ if and
only if $P_{\alpha}$ is not $t$-invertible for each $\alpha$.
\end{thm}
%%%%%%%%%%%%%%%%%%%%%%%%%%%%%%%%%%%%%%%%%%%%%%%%%%%%%%%%%%%
\begin{proof}
$(\Rightarrow)$ If $I^{-1}$ is a ring, by \cite[Proposition
2.1(2)]{20}, $(P_{\alpha})^{-1}$ is a ring for each $\alpha$. So, by
Lemma~\ref{sec:2.1}, $P_{\alpha}$ is not $t$-invertible for each $\alpha$.
Whence $\{P_{\alpha}\}\subseteq S_2\cup S_3$.\\

$(\Leftarrow)$ By \cite[Lemma 4.3]{20}, it is enough to prove that
$I^{-1}\subseteq(\bigcap R_{P_{\alpha}})\cap (\bigcap
R_{M_{\beta}})$ where $\{M_{\beta}\}$ is the set of all $t$-maximal
ideals of $R$ that do not contain $I$. Clearly $I^{-1}\subseteq \bigcap R_{M_{\beta}}$ (for if $x\in I^{-1}$ and $a\in I\backslash M_{\beta}$, then $x={xa\over a}\in R_{M_{\beta}}$). Now we show that $I^{-1}\subseteq \bigcap R_{P_{\alpha}}$. Let $P_{\alpha}$ be any minimal prime over $I$.
Since $P_{\alpha}$ is not $t$-invertible, $P_{\alpha}\in S_2\cup S_3$.
If $P_{\alpha}\in S_2$, set $J:=\bigcap_{\gamma\neq\alpha}P_{\gamma}$.
Then $I=J\cap P_{\alpha}$ and since $\bigcap P_{\alpha}$ is irredundant, $J\nsubseteq P_{\alpha}$.
But since $P_{\alpha}$ is a non-$t$-invertible $t$-maximal ideal of a $PVMD$ $R$, $(J+P_{\alpha})_{t}=R$ and $(P_{\alpha})^{-1}=R$.\\

\noindent{\bf Claim}. Let $R$ be a $PVMD$ and $A$ and $B$ nonzero ideals of $R$ such that $(A+B)_{t}=R$. Then $(A\cap B)_{t}=(AB)_{t}$.
Indeed, by \cite{25} it suffices to check that $(A\cap B)_{t}R_{M}=(AB)_{t}R_{M}$ for every $t$-maximal ideal $M$ of $R$. Let $M$ be a $t$-maximal ideal of $R$. Since $A$ and $B$ are $t$-comaximal, then either $A\nsubseteq M$ or $B\nsubseteq M$. Without loss of generality, we may assume that $A\nsubseteq M$. Hence, by \cite[Lemma 3.3]{23} $(A\cap B)_{t}R_{M}=(A\cap B)R_{M}=AR_{M}\cap BR_{M}=R_{M}\cap BR_{M}=BR_{M}=ABR_{M}=(AB)_{t}R_{M}$, as desired.\\
Now, by the claim $I=J\cap P_{\alpha}=(J\cap P_{\alpha})_{t}=(JP_{\alpha})_{t}$. So
$I^{-1}=(JP_{\alpha})^{-1}=(R:P_{\alpha}J)=((R:P_{\alpha}):J)=(R:J)=J^{-1}$. But since $J\nsubseteq P_{\alpha}$, $I^{-1}=J^{-1}\subseteq
R_{P_{\alpha}}$. Assume that $P_{\alpha}\in S_3$ and let $N$ be a $t$-maximal ideal of $R$
properly containing $P_{\alpha}$. Since $I$ is a radical ideal of
$R$, $IR_N=P_{\alpha}R_N$. Since $P_{\alpha}R_N$ is a
nonmaximal prime ideal of the valuation domain $R_{N}$, it is not
invertible. Hence $I^{-1}\subseteq (I^{-1})_{R\backslash N}\subseteq
(R_N:IR_N)=(R_N:P_{\alpha}R_N)=R_{P_{\alpha}}$ (\cite[Corollary 3.6]{22}), as desired.
\end{proof}
\bigskip
%%%%%%%%%%%%%%%%%%%%%%%%%%%%%%%%%%%%%%%%%%%%%%%%%%%%%%%%%%%%%%%%%%%%%%%%%%%%%%%%%%%%%%%%%%%%%%%%%%
The following example shows that the irredundancy condition in Theorem~\ref{sec:2.3} cannot be removed.
This example is a slight modification of \cite[Example 5.1]{20}, where the authors
constructed an example of a Bezout domain $R$ with a principal ideal $I$ (so $I^{-1}$ is not a ring) such that $P^{-1}$ is a ring for each minimal prime ideal $P$ of $I$. Our example is just an adjunct of an indeterminate $Y$ to the domain $R$ to get outside the Pr\"ufer situation but keeping us in the context of $PVMDs$.\\
%%%%%%%%%%%%%%%%%%%%%%%%%%%%%%%%%%%%%%%%%%%%%%%%%%%%%%%%%%%%%%%%%%%%%%%%%%%%%%%%%%%%%%%%%%%%%%%%%%%%%%%%%%%%%%%%%%%%%%%%%%
\begin{example}\label{sec:2.4}
Let $\mathbb{Q}$ be the filed of rational numbers and set $T=\Q[\{X^n:n\in \Q^+\}]$ and $J=(X-1)T$. By
(\cite[Example 5.1]{20}), $T$ is a Bezout domain, $J$ is a principal
radical ideal of $T$ (so $J^{-1}$ is not a ring) and $P^{-1}$ is a
ring for each minimal $P$ over $J$ in $T$. Also, by \cite[Theorem 3.8]{22},
the intersection of the minimal primes of $J$ is not an irredundant
intersection. Now let $R=T[Y]$, $I=J[Y]$. Clearly $R$ is a $PVMD$ (which is not Pr\"ufer),
and $I$ is a radical principal ideal of $R$ (so $I^{-1}=J^{-1}[Y]$ is not a ring).
Since $J=I\cap T\subseteq Q\cap T=P$, it is easy to check that every minimal prime ideal $Q$ of $R$ over $I$ is of the form $Q=P[Y]$, where
$P$ is a minimal prime ideal of $T$ over $J$. Hence $Q^{-1}=P^{-1}[Y]$ is a ring for each $Q$. Finally
$I=J[Y]=(\bigcap P)[Y]=\bigcap P[Y]$ is not an irredundant intersection.
\end{example}
\bigskip
%%%%%%%%%%%%%%%%%%%%%%%%%%%%%%%%%%%%%%%%%%%%%%%%%%%%%%%%%%%%%%%%%%%%%%%%%%%%%%%%%
Let $T$ be an overring of an integral domain $R$. According to \cite{6}, $T$ is said to be
$t$-linked over $R$ if for each finitely generated ideal $I$ of $R$ with
$I^{-1}=R$, we have $(IT)^{-1}=T$. Also we say that
$T$ is $t$-flat over $R$ if $T_M=R_P$ for each $t$-maximal ideal
$M$ of $T$, where $P=R\cap M$ (cf. \cite{26}). Finally, we say that $Spec_{t}(R)$ is Noetherian
if $R$ satisfies the a.c.c. condition on the radical $t$-ideals.\\
Our second main theorem generalizes Heinzer-Papick's theorem \cite[Theorem 2.5]{18}.
\bigskip
%%%%%%%%%%%%%%%%%%%%%%%%%%%%%%%%%%%%%%%%%%%%%%%%%%%%%%%%%%%%%%%%%%%%%%%%%%%%%%%%%
%%%%%%%%%%%%%%%%%%%%%%%%%%%%%%%%%%%%%%%%%%%%%%%%%%%%%%%%%%%%%%%%%%%%%%%%%%%%%%%%%

\begin{thm}\label{sec:2.5}
Let $R$ be a $PVMD$ with $Spec_t(R)$ Noetherian, and let $I$ be a $t$-ideal of $R$. Assume that $I^{-1}$ is a ring. Then the
following conditions are equivalent:
\begin{enumerate}
    \item $I^{-1}=(I:I)$;
    \item $I=\sqrt{I}$;
    \item The minimal prime ideals of $I$ in $(I:I)$ are all $t$-maximal ideals.
\end{enumerate}
\end{thm}
%%%%%%%%%%%%%%%%%%%%%%%%%%%%%%%%%%%%%%%%%%%%%%%%%%%%%%%%%%%%%%%%%%%%%%%

The proof of this theorem involves several lemmas of independent interest, some of them are $t$-analogues of well-known results.
%%%%%%%%%%%%%%%%%%%%%%%%%%%%%%%%%%%%%%%%%%%%%%%%%%%%%%%%%%%%%%%%%%%%%%%

\begin{lemma}\label{sec:2.6}
Let $T$ be a $t$-flat overring of a domain $R$. The following equivalent conditions hold:
\begin{enumerate}
    \item $I_t\subseteq (IT)_{t_1}$ for each $I\in F(R)$, where $t_1$ is the $t$-operation w.r. to $T$.
    \item If $J$ is a $t$-ideal of $T$ and $J\cap R\neq 0$, then $J\cap R$ is
    a $t$-ideal of $R$.
    \item $I_vT\subseteqq (IT)_{v_1}$ for each $I\in f(R)$, where $v_1$ is the $v$-operation w.r. to $T$.
    \item $(IT)_{v_1}=(I_vT)_{v_1}$ for each $I\in f(R)$.
    \item $(IT)_{t_1}=(I_tT)_{t_1}$ for each $I\in F(R)$.
    \item $(IT)_{v_1}=(I_tT)_{v_1}$ for each $I\in F(R)$.
\end{enumerate}
\end{lemma}
%%%%%%%%%%%%%%%%%%%%%%%%%%%%%%%%%%%%%%%%%%%%%%%%%%%
\begin{proof}
The six conditions are equivalent for an arbitrary overring $T$ of $R$
by \cite[Proposition 1.1 ]{3}. To prove \i, let $x\in I_t$.
Then there is a finitely generated ideal $J$ of $R$ such that
$J\subseteq I$ and $x(R:J)\subseteq R$. Now, let $N$ be a $t$-maximal ideal of $T$ and set $M=N\cap R$. Since $T$ is $t$-flat over $R$,
$T_N=R_M$. Since $J$ is finitely generated, $x(T:JT)T_N=x(T_N:JT_N)=x(R_M:JR_M)=x(R:J)R_M\subseteq R_M=T_N$. Hence $x(T:JT)\subseteq T$ and so $x\in (JT)_{v_1}\subseteq (IT)_{t_1}$, as desired.
\end{proof}
%%%%%%%%%%%%%%%%%%%%%%%%%%%%%%%%%%%%%%%%%%%%%%%%%%%%%%%%%%%%%%%%%%%
%%%%%%%%%%%%%%%%%%%%%%%%%%%%%%%%%%%%%%%%%%%%%%%%%%%%%%%%%%%%%%%%%%%
%%%%%%%%%%%%%%%%%%%%%%%%%%%%%%%%%%%%%%%%%%%%%%%%%%%%%%%%%%%%%%%%%%%
%%%%%%%%%%%%%%%%%%%%%%%%% Lemma 3.3 %%%%%%%%%%%%%%%%%%%%%%%%%%%%%%%
The next lemma is crucial and it is a generalization of \cite[Theorem 26.1]{16}.
We will often use it whenever we want to prove that an overring $T$ of a $PVMD$ domain $R$
is contained in $R_{Q}$ for some $t$-prime ideal $Q$ of $R$.\\
%%%%%%%%%%%%%%%%%%%%%%%%%%%%%%%%%%%%%%%%%%%%%%%%%%%%%%%%%%%%%%%%%%%%%%%%%%%%%%%%%%%%%%%%%%%%
%%%%%%%%%%%%%%%%%%%%%%%%%%%%%%%%%%%%%%%%%%%%%%%%%%%%%%%%%%%%%%%%%%%%%%%%%%%%%%%%%%%%%%%%%%%%
\begin{lemma}\label{sec:2.7}
Let $R$ be a $PVMD$ and $T$ a $t$-linked overring of $R$. Then:
\begin{enumerate}
    \item If $M$ is a $t$-prime ideal of $T$, then $T_M=R_P$ and
    $M=PR_P\cap T$, where $P=M\cap R$.
    \item If P is a nonzero $t$-prime ideal of $R$, then
    $(PT)_{t_1}\neq T $ if and only if $R_P\supseteq T$, where
    $t_1$ is the $t$-operation w.r. to $T$.
    \item If $J$ is a $t$-ideal of $T$ and $I=J\cap R$, then
    $J=(IT)_{t_1}$.
    \item $\{(PT)_{t_1}\}_{P\in \Delta}$ is the set of all $t$-prime
    ideals of $T$, where $\Delta=\{P\in Spec_t(R):(PT)_{t_1}\neq
    T\}$.
\end{enumerate}
\end{lemma}
%%%%%%%%%%%%%%%%%%%%%%%%%%%%%%%%%%%%%%%%%%%%%%
\begin{proof}

    \i Since $T$ is a $t$-linked overring of a PVMD $R$, $T$
    is a $t$-flat overring of $R$ (\cite[Proposition 2.10]{26}).
    Hence $R_P=T_M$ where $P=M\cap R$ (\cite[Theorem 2.6]{8}). Therefore $M=MT_M\cap T=PR_P\cap T$.\\

    \ii If $(PT)_{t_1}\subsetneqq T$, then there is a
    $t$-maximal ideal $M$ of $T$ such that $M\supseteq (PT)_{t_1}$.
    Since $M\cap R\supseteq (PT)_{t_1}\cap R\supseteq PT\cap
    R\supseteq P$, $R_P\supseteq R_{M\cap R}=T_M\supseteq T$, as desired.\\
    Conversely, if $R_P\supseteq T\supseteq R$, then
    $T_{R\setminus P}=R_P$. Hence $R_P$ is $t$-linked over $T$. So, by Lemma~\ref{sec:2.6}, $(PT)_{t_1}\subseteq (PR_P)_{t_2}=PR_P\subsetneqq
    R_P$ (here $t_2$ is the $t$-operation w.r. to $R_P$ and it is trivial since $R_P$ is valuation). Since $T_{R\backslash P}=R_P$ is a valuation overring of a PVMD $T$, $J_{t_1}T_{R\backslash P}=JT_{R\backslash P}$ for each ideal $J$ of $T$. If $(PT)_{t_1}=T$, then $R_P=T_{R\backslash P}=(PT)_{t_1}T_{R\backslash P}=PT_{R\backslash P}=PR_{P}$, a contradiction. Therefore $(PT)_{t_1}\subsetneqq T$.\\

    \iii Clearly $(IT)_{t_1}\subseteq J$. It suffices to show that
    $J\subseteq (IT)_{t_1}$. Let $\{M_{\alpha}\}$ be the set of all
    $t$-maximal ideals of $T$. Since $T$ is a $t$-linked overring of $R$, $T$ is a PVMD. Hence $J=\bigcap
    JT_{M _{\alpha}}$. Set $P_{\alpha}=M_{\alpha}\cap R$ for each
    $\alpha$ and let $x\in JR_{M_\alpha}=JR_{P_{\alpha}}$. Then $x=\frac{a}{t}$, where $a\in J$ and $t\in
    R\setminus P_{\alpha}$. Since $J\subseteq T\subseteq
    T_{M_{\alpha}}=R_{P_{\alpha}}$, then $a=\frac{b}{s}$, where $b\in R$ and $s\in
    R\setminus P_{\alpha}$. Hence $b=as\in J\cap R=I$. So
    $x=\frac{b}{st}\in IR_{P_{\alpha}}\subseteq(IT)R_{P_{\alpha}}=(IT)T_{M_{\alpha}}$.
    Therefore $J\subseteq (IT)_{t_1}$, as desired.\\

    \iv By \iii, each $t$-prime ideal of $T$ is of the form $(PT)_{t_1}$
    for some $P\in\Delta$. Conversely, if $P\in\Delta$, then $P_tR_P=PR_P$
    is a $t$-prime ideal of $R_P$ (\cite[Lemma 3.3]{23} and $R_P$ is a valuation domain) and
    $T\subseteq R_P$ (by part(ii)). So $R_P=T_{R\setminus P}$ and then $R_P$ is $t$-linked over $T$.
    Hence $PR_P\cap T$ is a $t$-prime ideal of $T$ ( Lemma~\ref{sec:2.6}) and $PR_P\cap T=(((PR_P\cap T)\cap R)T)_{t_1}=(PT)_{t_1}$
    by \iii.
\end{proof}
\bigskip
%%%%%%%%%%%%%%%%%%%%%%%%%%%%%%%%%%%%%%%%%%%%%%%%%%%%%%%%%%%%%%%%%%%
%%%%%%%%%%%%%%%%%%%%%%%%%%%%%%%% Lemma 3.4 %%%%%%%%%%%%%%%%%%%%%%%%
The next lemma is a generalization of \cite[Lemma 2.4]{18} and it relates the fact $I^{-1}$ not being a ring to a kind of
``separation property" for a minimal prime ideal over a $t$-ideal of a $PVMD$.\\
%%%%%%%%%%%%%%%%%%%%%%%%%%%%%%%%%%%%%%%%%%%%%%%%%%%%%%%%%%%%%%%%%%%%%%%%%%%%%%%%%%%%%%
%%%%%%%%%%%%%%%%%%%%%%%%%%%%%%%%%%%%%%%%%%%%%%%%%%%%%%%%%%%%%%%%%%%%%%%%%%%%%%%%%%%%%%
\begin{lemma}\label{sec:2.8}
Let $R$ be a $PVMD$, $I$ a $t$-ideal of $R$ and $P$ a minimal
prime ideal over $I$ in $R$. If there is a finitely generated ideal $J$ of $R$
such that $I\subseteq J\subseteq P$, then $I^{-1}$ is not a ring.
\end{lemma}
%%%%%%%%%%%%%%%%%%%%%%%%%%%%%%%%%%%%%%%%%%%%%%%%%%%
\begin{proof}
By the way of contradiction, assume that $I^{-1}$ is a ring. Then by
\cite[Theorem 4.5]{20} and \cite[Theorem 2.22]{25}, $I^{-1}\subseteq
R_P$ and $I^{-1}$ is a $t$-linked overring of $R$. So $R_P$ is $t$-linked
over $I^{-1}$. Since $J^{-1}\subseteq I^{-1}$, $R=(JJ^{-1})_t\subseteq(JI^{-1})_{t_1}$
where $t_1$ is the $t$-operation w.r. to $I^{-1}$ (Lemma~\ref{sec:2.6}).
Also by Lemma~\ref{sec:2.6}, $(PI^{-1})_{t_1}\subseteq(PR_P)_{t_2}=PR_P$ (where $t_2$ is the
$t$-operation w.r. to $R_P$, so it is trivial). Therefore
$1\in R=(JJ^{-1})_t\subseteq (JI^{-1})_{t_1}\subseteq
(PI^{-1})_{t_1}\subseteq PR_P$, which is a contradiction.
\end{proof}
\bigskip
%%%%%%%%%%%%%%%%%%%%%%%%%%%%%%%%%%%%%%%%%%%%%%%%%%%%%%%%%%%%%%%%%%%
%%%%%%%%%%%%%%%%%%%%%%%%%%%%%% Lemma 3.5 %%%%%%%%%%%%%%%%%%%%%%%%%%
\begin{lemma}{\rm (\cite[Lemma 2.6]{24})}\label{sec:2.9}
Let $R$ be a $PVMD$ and $I$ a $t$-ideal of $R$. Then $I$ is a
$t$-ideal of $(I:I)$.
\end{lemma}
%%%%%%%%%%%%%%%%%%%%%%%%%%%%%%%%%%%%%%%%%%%%%%%%%%%%%%%%%%%%%%%%%%%
%%%%%%%%%%%%%%%%%%%%%%%%%%%%%%%%%%%%%%%%%%%%%%%%%%%%%%%%%%%%%%%%%%%
%%%%%%%%%%%%%%%%%%%%%%%%%%%%%%%%%%%%%%%%%%%%%%%%%%%%%%%%%%%%%%%%%%%
%%%%%%%%%%%%%%%%%%%%%%%% Lemma 3.6 %%%%%%%%%%%%%%%%%%%%%%%%%%%%%%%%
\begin{lemma}{\rm (\cite[Lemma 3.7)]{8}}\label{sec:2.10}
Let $R$ be an integral domain. The following conditions are equivalent.\\
\i Each $t$-prime ideal is the radical of a $v$-finite ideal.\\
\ii Each radical $t$-ideal is the radical of a $v$-finite ideal.\\
\iii $Spec_t(R)$ is Noetherian.
\end{lemma}
\bigskip
%%%%%%%%%%%%%%%%%%%%%%%%%%%%%%%%%%%%%%%%%%%%%%%%%%%%%%%%%%%%%%%%%%%
%%%%%%%%%%%%%%%%%%%%%%%%%%% Proposition 3.7 %%%%%%%%%%%%%%%%%%%%%%%
%%%%%%%%%%%%%%%%%%%%%%%%%%%%%%%%%%%%%%%%%%%%%%%%%%%%%%%%%%%%%%%%%%
\noindent{\bf Proof of Theorem 2.5}
$(ii)\Rightarrow(i)$  Follows from \cite[Proposition 3.3]{1} without any more conditions.\\

$(i)\Rightarrow(ii)$ Deny, assume that $I\subsetneqq\sqrt{I}$. Then there is a $t$-maximal ideal $M$ of $R$ such that $IR_M$ is not
a radical ideal. Moreover, there is a prime ideal $P$ contained in
$M$ and minimal over $I$ with $IR_M\subsetneqq PR_M$ and
$\sqrt{IR_M}=PR_M$. Note that $P$ is a $t$-prime ideal of $R$ (as a minimal prime over a $t$-ideal).
\begin{claim} $IR_P=PR_P$.\end{claim}
Deny. Let $b\in P$ such that $IR_P\subsetneqq bR_P\subseteq PR_P$.
Since $Spec_t(R)$ is Noetherian, $P=\sqrt{(a_1,...,a_r)_v}$ for some
$a_1,..., a_r\in P$. Set $J:=(a_1,...,a_r,b)$. Note that
$P=\sqrt{J_v}$ ($P=\sqrt{(a_1,...,a_r)_v}\subseteq\sqrt{(a_1,..., a_r,b)_v}\subseteq
P$). Now, we prove that $I\subseteq J\subseteq P$, which contradicts the assumption that
$I^{-1}$ is a ring by Lemma~\ref{sec:2.8}. Let $N$ be a $t$-maximal
ideal of $R$. If $P\nsubseteq N$, then
$R_N=PR_N=\sqrt{J_vR_N}=\sqrt{J_tR_N}=\sqrt{JR_N}$ (the last
equality holds since $N$ is $t$-prime, \cite[Lemma
3.3]{23}). Hence $JR_N=R_N\supseteq IR_N$. Assume that $P\subseteq N$. Then
$PR_P=PR_N$ since $R_P$ is an overring of the valuation domain $R_N$. Since $IR_P\subsetneqq bR_P$, $b^{-1}I\subsetneqq R_P$
and so $b^{-1}I\subseteq PR_P=PR_N\subseteq R_N$. Hence
$IR_N\subseteq bR_N\subseteq JR_N$ as desired.\\
Now since $R_M$ is a valuation domain, $Z(R_M,IR_M)=QR_M$ for some
$t$-prime ideal $Q\subseteq M$. Since $R$ is a $PVMD$ and $P$ and $Q$ are $t$-primes under
the $t$-maximal ideal $M$, $Q$ and $P$ are comparable under inclusion. Moreover, let $x\in
PR_M\setminus IR_M$. Since $PR_M=PR_P=IR_P$ (Claim 1), there exists $y\in R\setminus P$ such that $yx\in
I$. Hence $y\in Z(R_M, IR_M)\cap R=Q$ and therefore $P\subsetneqq Q$.\\
\begin{claim} $(QI^{-1})_{t_1}=I^{-1}$.\end{claim}
Note that $I^{-1}=(I:I)$ is a subintersection of $R$ (\cite[Theorem 4.5]{20}) and so $I^{-1}$ is $t$-linked
over $R$ (\cite[Theorem 2.22]{25}). Since $Spec_t(R)$ is
Noetherian, $Q=\sqrt{A_v}$ for some finitely generated ideal $A$ of $R$. Say
$A=\displaystyle\sum^{n=m}_{n=1}b_nR$. Since $P\subsetneqq Q$, $P\subsetneqq A_v$. Indeed, let $N$ be a $t$-maximal ideal of $R$.
If $Q\nsubseteq N$, then $PR_N\subseteq R_N=QR_N=AR_N$.  If $Q\subseteq N$, then $AR_N$ and $PR_N$ are comparable as ideals of the valuation domain $R_N$. But if $AR_N\subseteq PR_N$, then $QR_N=\sqrt{A_vR_N}=\sqrt{A_tR_N}=\sqrt{AR_N}\subseteq PR_N$ and so $Q\subseteq P$, which is absurd. Hence
$PR_N\subsetneqq AR_N$ and therefore $P\subsetneqq A_t=A_v$. Now since $I\subseteq P\subseteq A_v$,
$A^{-1}\subseteq I^{-1}$. So $1\in R=(AA^{-1})_t\subseteq (AI^{-1})_t\subseteq
(AI^{-1})_{t_1}\subseteq (QI^{-1})_{t_1}$ (Lemma~\ref{sec:2.6}). Hence
$(QI^{-1})_{t_1}=I^{-1}$, as desired.\\
Finally, by Lemma~\ref{sec:2.7}, $I^{-1}\nsubseteq R_Q$. On the other hand $(I:I)\subseteq
(I:I)R_M\subseteq(IR_M:IR_M)=(R_M)_{QR_M}=R_Q$ by \cite[Lemma 3.1.9]{13}, which is absurd. It follows that $I$ is a radical ideal of $R$.\\

$(iii)\Rightarrow(ii)$ Assume that all minimal prime ideals of $I$ in $(I:I)$ are
$t$-maximal ideals. If $I\subsetneqq\sqrt{I}$, as in the proof of $(i)\Rightarrow(ii)$, there exist
two $t$-prime ideals $P$ and $Q$ of $R$ such that $I\subseteq
P\subsetneqq Q$ and $(I:I)\subseteq R_Q$. Then $(I:I)_{R\setminus Q}=R_Q$ and so
$R_Q$ is $t$-linked over $(I:I)$. Hence $QR_Q\cap (I:I)$ and $PR_Q\cap
(I:I)$ are $t$-prime ideals of $(I:I)$ with $I\subseteq
PR_Q\cap (I:I)\subsetneqq QR_Q\cap (I:I)$ which is absurd.\\

$(i)\Rightarrow(iii)$ Assume that $I^{-1}=(I:I)$ and let $P$ be a minimal
prime of $(I:I)$ over $I$. By Lemma~\ref{sec:2.9}, $I$ is a $t$-ideal of $(I:I)$ and so $P$ is a $t$-prime ideal of $(I:I)$ (as a minimal prime over a $t$-ideal). Now by a way of contradiction, assume that there is a
$t$-prime ideal $Q$ of $(I:I)$ such that $P\subsetneqq Q$. Since
$(I:I)$ is  a $t$-linked overring of $R$, $P=(P'(I:I))_{t_1}$
and $Q=(Q'(I:I))_{t_1}$ for some $t$-prime ideals $P'$ and $Q'$ of
$R$ with $I\subseteq P'\subsetneqq Q'$ (Lemma~\ref{sec:2.7}(iv)). Set $Q'=\sqrt{A}$ for some finitely generated ideal $A$ of $R$. As in the
proof of Claim 2, $I\subseteq P'\subseteq A_{t}$. So $A^{-1}\subseteq I^{-1}=(I:I)$. Hence $1\in R=(AA^{-1})_t\subseteq (A(I:I))_{t_1}\subseteq (Q'(I:I))_{t_1}=Q$, which is absurd. It follows that $P$ is a $t$-maximal ideal of $(I:I)$, completing the proof. $\square$\\

\medskip
%%%%%%%%%%%%%%%%%%%%%%%%%%%%%%%%%%%%%%%%%%%%%%%%%%%%%%%%%%%%%%%%%%
%%%%%%%%%%%%%%%%%%%%%%%%%%%%% Proposition 3.8 %%%%%%%%%%%%%%%%%%%%
%%%%%%%%%%%%%%%%%%%%%%%%%%%%%%%%%%%%%%%%%%%%%%%%%%%%%%%%%%%%%%%%%%
The next two results deal with the duals of primary $t$-ideals in a $PVMD$.
%%%%%%%%%%%%%%%%%%%%%%%%%%%%%%%%%%%%%%%%%%%%%%%%%%%%%%%%%%%%%%%%%%%%%%%%%%%

\begin{proposition}{\rm (cf. \cite[Lemma 4.4]{12})}\label{sec:2.11}
Let $R$ be a $PVMD$ and $I$ a primary $t$-ideal of $R$. If
$I^{-1}$ is a ring, then $I^{-1}=(I:I)$.
\end{proposition}
%%%%%%%%%%%%%%%%%%%%%%%%%%%%%%%%%%%%%%%%%%%%%%
\begin{proof}
Deny, assume that there is $x\in I^{-1}\backslash (I:I)$. Since $I$
is a $t$-ideal of $R$, there is $a\in I$ and a $t$-maximal ideal
$M$ of $R$ such that $xa\not\in IR_M$. Since $I^{-1}=(\bigcap
R_{P_{\alpha}})\cap(\bigcap R_{M_{\beta}})$ (\cite[Theorem 4.5]{20}),
$x\in R_{M_{\beta}}$ for each $\beta$ and hence $I\subseteq M$.
Therefore there is a minimal prime $I\subseteq P_{\alpha}\subseteq M$. Thus
$x\in R_{P_{\alpha}}$. Write $x=\frac{b}{s}$ where $b\in R$ and
$s\in R\backslash P_{\alpha}$. If $t=\frac{s}{a}\in R_M$, then
$s=ta\in PR_M\cap R=P$, which is a contradiction. If $\frac{a}{s}\in R_M$,
since $I$ is a primary ideal of $R$, $ax=a\frac{b}{s}=b\frac{a}{s}\in
IR_{P_{\alpha}}\cap R_M=IR_M$, which is a contradiction too. It follows that $I^{-1}=(I:I)$.
\end{proof}
\bigskip
%%%%%%%%%%%%%%%%%%%%%%%%%%%%%%%%%%%%%%%%%%%%%%%%%%%%%%%%%%%%%%%%%%%
%%%%%%%%%%%%%%%%%%%%%%%%%%% Corollary 3.9 %%%%%%%%%%%%%%%%%%%%%%%%%
\begin{corollary}{\rm (cf. \cite[Proposition 3.1.14]{13})}\label{sec:2.12}
Let $R$ be a $PVMD$ with $Spec_t(R)$ Noetherian and $I$ a
$t$-ideal of $R$. If $I$ is a primary ideal which is not prime, then $I^{-1}$ is
not a ring .
\end{corollary}

\begin{proof}
Deny, assume that $I^{-1}$ is a ring. Then $I^{-1}=(I:I)$ by Proposition~\ref{sec:2.11}.
Therefore $I$ is a radical ideal (and so prime) by Theorem~\ref{sec:2.5}, which is absurd.
\end{proof}
\bigskip
%%%%%%%%%%%%%%%%%%%%%%%%%%%%%%%%%%%%%%%%%%%%%%%%%%%%%%%%%%%%%%%%%%%%%%%%%%%%%%%%%%%%%%%%%%%%%%%%%%%%%%
%%%%%%%%%%%%%%%%%%%%%%%%%%%%%%%%%%%%%%%%%%%%%%%%%%%%%%%%%%%%%%%%%%%%%%%%%%%%%%%%%%%%%%%%%%%%%%%%%%%%%%
According to \cite[Section 27]{16}, a Pr\"uer domain $R$ is called a QR-domain
if each overring of $R$ is a quotient ring of $R$. In \cite{7} the
authors defined $t$QR-domains as $PVMDs$ $R$ such that each $t$-linked
overring of $R$ is a quotient ring of $R$ and they characterized
$t$QR-domains as follows: ``{\em Let $R$ be a $PVMD$. Then $R$ is a $t$QR-domain if and only if for each
f.g. ideal $A$ of $R$, there is $n\geq 1$ and $b\in R$ such that
$A^n\subseteq bR\subseteq A_v$}" \cite[Theorem 1.3]{7}.\\
We close this section with a third main theorem. It generalizes a well-known results by Fontana et al. \cite[Corrollary 4.4 and Theorem 4.11]{14} and gives a description of $(I:I)$ for a $t$-ideal $I$ in a $PVMD$ that is a $t$QR-domain.
%%%%%%%%%%%%%%%%%%%%%%%%%%%%%%%%%%%%%%%%%%%%%%%%%%%%%%%%%%%%%%%%%%%

\begin{thm}\label{sec:2.13}
Let $I$ be a $t$-ideal of a $PVMD$ $R$, $\{Q_{\alpha}\}$
be the set of all maximal prime ideals of $Z(R,I)$ and
$\{M_{\beta}\}$ be the set of $t$-maximal ideals of $R$ that do not
contain $I$. Then:
\begin{enumerate}
    \item $(I:I)\supseteq(\bigcap R_{Q_{\alpha}})\cap(\bigcap
    R_{M_{\beta}})$;
    \item If $R$ is a $tQR$-domain then the equality holds.
\end{enumerate}
\end{thm}
%%%%%%%%%%%%%%%%%%%%%%%%%%%%%%%%%%%%%%%%%%%%%%%%%%%%%%%%%%%%%%%%%%%%%%%
%%%%%%%%%%%%%%%%%%%%%%%%%%%%%%%%%%%%%%%%%%%%%%%%%%%%%%%%%%%%%%%%%%%%%%%
Before proving this theorem, we need the following lemma.
\bigskip

%%%%%%%%%%%%%%%%%%%%%%%% Proposition 4.1 %%%%%%%%%%%%%%%%%%%%%%%%%%
\begin{lemma}\label{sec:2.14}
Let $I$ be a $t$-ideal of a $PVMD$ $R$ and let $\{Q_{\alpha}\}$
be the set of all prime ideals of $Z(R,I)$. Then
$Q_{\alpha}$ is a $t$-prime ideal for each $\alpha$.
\end{lemma}
%%%%%%%%%%%%%%%%%%%%%%%%%%%%%%%%%%%%%%%%%%%%%%%%%%%%
\begin{proof}
First we claim that $Z(R,I)=\displaystyle\bigcup_{M\in M_t(R,I)}Z(R_M,IR_M)\cap R$.
Indeed, let $x\in Z(R,I)$. Then there is $a\in R\backslash I$ such that $ax\in
I$. Since $I$ is a $t$-ideal, there is a $t$-maximal ideal $I\subseteq M$ of
$R$ such that $a\in R_M\backslash IR_M$ and $ax\in IR_M$. Hence
$x\in Z(R_M,IR_M)\cap R$. Conversely, let $M\in Max_t(R,I)$. Since $R_M$ is a
valuation domain, there is a $t$-prime ideal $Q\subseteq M$
such that $Z(R_M,IR_M)=QR_M$. Now we prove that $Q\subseteq Z(R,I)$.
Let $z\in Q$. Then $z\in QR_M$ and hence there is $\frac{c}{t}\in
R_M\backslash IR_M$ such that $\frac{zc}{t}\in IR_M$ with $c\in
R\backslash I$ and $t\in R\backslash M$. This implies that $szc\in
I$ for some $s\in R\backslash M$. If $cs\in I$, then $c=\frac{i}{s}\in IR_M$. Thus $\frac{c}{t}\in IR_M$, a contradiction. Then $cs\not\in I$ and then $z\in Z(R,I)$. Therefore $Z(R_M,IR_M)\cap R=QR_M\cap R=Q\subseteq Z(R,I)$. Finally, $Q$'s are $t$-prime ideals of $R$ (\cite[Corollary 2.47]{25}).
\end{proof}
%%%%%%%%%%%%%%%%%%%%%%%%%%%%%%%%%%%%%%%%%%%%%%%%%%%%%%%%%%%%%%%%%%%
%%%%%%%%%%%%%%%%%%%%%%%% Proposition 4.1 %%%%%%%%%%%%%%%%%%%%%%%%%%
%%%%%%%%%%%%%%%%%%%%%%%%%%%%%%%%%%%%%%%%%%%%
\noindent{\bf Proof of Theorem 2.13}.
\i Let $u\in (\bigcap R_{Q_{\alpha}})\cap(\bigcap
    R_{M_{\beta}})$ and $a\in I$. It is enough to prove that $ua\in I$. Since $u\in \bigcap R_{M_{\beta}}$,
    it suffices to show that $ua\in R_{N_{\gamma}}$ for each $\gamma$, where $\{N_{\gamma}\}$
    be the set of $t$-maximal ideals of $R$ containing
    $I$. By \cite[Corollary 4.6]{16}, $\bigcap R_{Q_{\alpha}}=R_{R\setminus\cup
    Q_{\alpha}}$. Write $u=\frac{r}{s}$, where $r\in R$ and $s\in R\setminus\cup
    Q_{\alpha}$. Fix $\gamma$ and choose $\alpha_1$ such that
    $I\subseteq Q_{\alpha_1}\subseteq N_{\gamma}$. We claim that
    $\frac{a}{s}\in R_{N_{\gamma}}$. For if not, then
    $\frac{s}{a}=t\in R_{N_{\gamma}}$ and thus $s=at\in
    Q_{\alpha_1}R_{N_{\gamma}}\cap R=Q_{\alpha_1}$, a contradiction.
    If $ua\not\in IR_{N_{\gamma}}$, then
    $ua=r(\frac{a}{s})=\frac{c}{b}$, where $c\in R\setminus I$ and
    $b\in R\setminus N_{\gamma}$. Hence $sc=rab\in I$. Thus
    $s\in\cup Q_{\alpha}$, a contradiction. Therefore $ua\in IR_{N_{\gamma}}$, as desired.\\

\ii Set $T:=(I:I)$.  Clearly $T\subseteq \bigcap R_{M_{\beta}}$. Now we
    will prove that $T\subseteq \bigcap R_{Q_{\alpha}}$. By Lemma~\ref{sec:2.7}(ii), it
    suffices to show that $(Q_{\alpha}T)_{t_1}\neq T$ for each
    $\alpha$. Since $R$ is a $PVMD$ and $I$ is a $t$-ideal, $T$ is $t$-linked over $R$ .
    Hence $T=R_S$ for some multiplicative closed set $S$ of $R$ since $R$ is a $tQR$-domain.
    By the way of contradiction, assume that $(QT)_{t_1}=T$ where $Q=Q_{\alpha}$ for some $\alpha$.
    Then there exists a finitely generated ideal $B$ such that $B_{v_1}=T$ and $B\subseteq QT$.
    Say $B=\displaystyle\sum^{i=r}_{i=1}a_nT$ with $a_i\in QT$ and write $a_i=\displaystyle\sum^{s=m_i}_{s=1}q_{is}t_{is}$ with $q_{is}\in Q$ and
    $t_{is}\in T$ for each $i=1, \dots, n$ and $s=1, \dots, m_i$. Now let
    $A$ be the finitely generated ideal of $R$ generated by all $q_{is}'s$. Then $A\subseteq Q$ and $B\subseteq AT$. Hence $T=B_{v_1}\subseteq (AT)_{v_1}\subseteq (A_vT)_{v_1}\subseteq T$ and therefore $(AT)_{v_1}=(A_vT)_{v_1}=T$. Since $R$ is a $tQR$-domain and $T$ is $t$-linked over $R$, by \cite[Proposition 2.17]{6}, $A_vT=T$. But since $A_v=A_t\subseteq Q$ (here $Q$ is a $t$-prime ideal by Lemma~\ref{sec:2.14}), $QT=T$.
    Hence $1=\displaystyle\sum^{i=n}_{i=1}q_ia_i$ where $q_i\in Q$ and $a_i\in T$. Set $J=\displaystyle\sum_{i=1}^{i=n}q_iR$. Clearly $JT=T$ and by induction $J^{s}T=T$ for all positive integer $s$. Since $R$ is a $tQR$-domain, there is a positive integer $N$ and $d\in R$ such
    that $J^{N}\subseteq dR\subseteq J_v=J_t\subseteq Q$. Since $J^{N}T=T$, then $1=\displaystyle\sum_{i=1}^{i=s}\lambda_iy_i$ where $\lambda_i\in J^{N}$ and $y_i\in T$, and since $J^{N}\subseteq dR$, there exists $\mu_i\in R$ such that $\lambda_i=d\mu_i$ for each $i$. Now, since $d\in Q\subseteq Z(R, I)$, there exists $r\in R\setminus I$ such that $rd\in I$. Hence $r=\displaystyle\sum_{i=1}^{i=s}r\lambda_iy_i=\displaystyle\sum_{i=1}^{i=s}rdy_i\mu_i\in IT=I$, a contradiction. Hence $(QT)_{t_1}\subsetneqq T$ and by Lemma~\ref{sec:2.7}, $T\subseteq R_{Q}$, completing the proof. $\square$

%%%%%%%%%%%%%%%%%%%%%%%%%%%%%%%%%%%%%%%%%%%%%%%%%%%%%%%%%%%%%%%%%%%%%%%%%%%%%%%%%%%
\section{Ideal Transform overrings of a $PVMD$}\label{sec:3}
%%%%%%%%%%%%%%%%%%%%%%%%%%%%%%%%%%%%%%%%%%%%%%%%%%%%%%%%%%%%%%%%%%%
%%%%%%%%%%%%%%%%%%%  Proposition 2.1 %%%%%%%%%%%%%%%%%%%%%%%%%%%%%%
%%%%%%%%%%%%%%%%%%%%%%%%%%%%%%%%%%%%%%%%%%%%%%%%%%%%%%%%%%%%%%%%%%%
We start this section with the following theorem which is a generalization of \cite[Theorem 3.2.5]{13}. As the proof
is similar to that one of \cite[Theorem 3.2.5]{13} simply by replacing maximal ideals by $t$-maximal ideals, we remove it here.
%%%%%%%%%%%%%%%%%%%%%%%%%%%%%%%%%%%%%%%%%%%%%%%%%%%%%%%%%%%%%%%%%%%%%%%%%%%%%%%%%%%%
\begin{thm}
Let $R$ be a $PVMD$, $I$ a $t$-ideal of $R$, $\{P_{\alpha}\}$ the
set of minimal prime ideals of $I$, and $\{M_{\beta}\}$ the set of
$t$-maximal ideals of $R$ that do not contain $I$. Then:
\begin{enumerate}
    \item $T(I)\subseteq(\bigcap R_{Q_{\alpha}})\cap(\bigcap R_{M_{\beta}})$,
    where $Q_{\beta}$ is the unique prime ideal determined by $\bigcap^{\infty}_{n=1} I^nR_{P_{\alpha}}$;
    \item The equality holds, if $I$ has a finitely many minimal primes.
\end{enumerate}
\end{thm}
\bigskip
%%%%%%%%%%%%%%%%%%%%%%%%%%%%%%%%%%%%%%%%%%%%%%%%%%%%%%%%%%%%%%%%%%%%%%%%%%%%%%%%%%%%%%%
Our next theorem generalizes \cite[Theorem 3.3.7]{13} to $PVMDs$.
%%%%%%%%%%%%%%%%%%%%%%%%%%%%%%%%%%%%%%%%%%%%%%%%%%%%%%%%%%%%%%%%%%%%%%%%%%%%%%%%%%%%%%%
%%%%%%%%%%%%%%%%%%%%%%%%%%%%%%%%%%%%%%%%%%%%%%%%%%%%%%%%%%%%%%%%%%%%%%%%%%%%%%%%%%%
\begin{thm}\label{sec:3.2}
Let $P$ be a non-$t$-invertible $t$-prime ideal of a $PVMD$ $R$. Then
there is no proper overring of $R$ between $P^{-1}$ and $\Omega(P)$.
\end{thm}
%%%%%%%%%%%%%%%%%%%%%%%%% Lemma 5.1 %%%%%%%%%%%%%%%%%%%%%%%%%%%%%%%%%%%%%%%%%%%%%%%%%%%%%%%%%%%%%%%%%%%%%%%%%%%%%%%%%%%%%%%%%%%%%%%%%%%%%%%%%%%%%%%%%%
The proof of this theorem involves the following lemmas.
\bigskip
%%%%%%%%%%%%%%%%%%%%%%%%%%%%%%%%%%%%%%%%%%%%%%%%%%%%%%%%%%%%%%%%%%%%%%%%%%%%%%%%%%%%%%%%%%%%%%%%%%%%%%%%%%%%%%%%%%%%%%%%%%%%%%%%%%%%%%%%%%%%%%%%%%%%%%%

\begin{lemma}\label{sec:3.3}
Let $R$ be a $PVMD$, $I$ a $t$-ideal of $R$ and let $T$ be a
$t$-linked overring of $R$ contained in $\Omega(I)$.
Then there is a one-to-one correspondence between the sets $S_1=\{P\in
Spec_t(R):P\nsupseteq I\}$ and $S_2=\{Q\in Spec_t(T):Q\nsupseteq IT\}$.
\end{lemma}

\begin{proof} Define $\Psi:S_1\rightarrow S_2$ by $\Psi(P)=PR_P\cap T=Q$ for each $P\in S_1$. Then $\Psi$ is well-defined. Indeed, let $P\in S_1$. Since $T\subseteq \Omega(I)$, $T\subseteq R_P$. So $T_{R\setminus P}=R_P$ and then $R_P$ is a $t$-linked overring of $T$. Hence $PR_P\cap T$ is a $t$-prime of $T$. Also, if $x\in I\setminus P$, then $x\in IT\setminus Q$ and the injectivity of $\Psi$ is clear.\\
Now, let $Q\in S_2$ and set $P:=R\cap Q$. Then $P\nsupseteq I$, and since $R_P=T_Q$, $PR_P=QT_Q$. Hence $\Psi(P)=PR_P\cap T=QT_Q\cap T=Q$.
\end{proof}
\bigskip
%%%%%%%%%%%%%%%%%%%%%%%%%%%%%%%%%%%%%%%%%%%%%%%%%%%%%%%%%%%%%%%%%%%
%%%%%%%%%%%%%%%%%%%%%%%%%% Lemma 5.2 %%%%%%%%%%%%%%%%%%%%%%%%%%%%%%
\begin{lemma}\label{sec:3.4}
Under the same notation as Lemma~\ref{sec:3.3}, if $(IT)_{t_1}=T$ , then $T=\Omega(I)$.
\end{lemma}

\begin{proof}
Assume that $(IT)_{t_1}=T$. Then $IT$ is not contained in any
$t$-prime ideal of $T$. Since $R$ is a $PVMD$ and $T$ is a $t$-linked
overring of $R$, $T$ is a $PVMD$. By Lemma~\ref{sec:3.3}, $T=\displaystyle\bigcap_{Q\in
Spec_t(T)} T_Q=\displaystyle\bigcap_{P\in Spec_t(R),P\nsupseteq I}
R_P\supseteq\Omega(I)$. Hence $T=\Omega(I)$.
\end{proof}
\bigskip
%%%%%%%%%%%%%%%%%%%%%%%%%%%%%%%%%%%%%%%%%%%%%%%%%%%%%%%%%%%%%%%%%%%
%%%%%%%%%%%%%%%%%%%%%%%% Corollary 5.3 %%%%%%%%%%%%%%%%%%%%%%%%%%%%
%%%%%%%%%%%%%%%%%%%%%%%%%%%%%%%%%%%%%%%%%%%%%%%%%%%%%%%%%%%%%%%%%%
%By corollary~\ref{sec:3.6}, $(IT(I))_{t_1}=T(I)$ implies that
%$(I\Omega(I))_{t_1}=\Omega(I)$. The converse is not true in general.
%For example: Let $(V,M)$ be a nondiscrete rank one valuation domain.
%Then $M=M^2$, and hence $T(M)=M^{-1}=V$ and $\Omega(M)=qf(V)$.
%%%%%%%%%%%%%%%%%%%%%%%%%%%%%%%%%%%%%%%%%%%%%%%%%%%%%%%%%%%%%%%%%%%
%%%%%%%%%%%%%%%%%%%%%%%% Lemma 5.4 %%%%%%%%%%%%%%%%%%%%%%%%%%%%%%%%

%%%%%%%%%%%%%%%%%%%%%%%%%%%%%%%%%%%%%%%%%%%%%%%%%%%%%
\noindent{\bf Proof of Theorem 3.2}.
Let $T$ be an overring of $R$ such that $P^{-1}\subsetneqq
T\subseteq\Omega(P)$ and let $\{M_{\beta}\}$ be the set of all $t$-maximal
ideals of $R$ that do not contain $P$. By \cite[Theorem 3.2.2]{13},
$T\subseteq\Omega(P)\subseteq\displaystyle\bigcap R_{M_{\beta}}$. If
$(PT)_{t_1}\neq T$, then $T\subseteq R_P$ (Lemma~\ref{sec:2.7}(ii)). So
$P^{-1}\subsetneqq T\subseteq R_P\cap(\bigcap R_{M_{\beta}})=P^{-1}$
(\cite[Proposition 1.2]{19}), which is a contradiction. Hence $(PT)_{t_1}=T$, and so $T=\Omega(P)$ by Lemma~\ref{sec:3.4}.$\square$
\bigskip
%%%%%%%%%%%%%%%%%%%%%%%%%%%%%%%%%%%%%%%%%%%%%%%%%%%%%%%%%%%%%%%%%%%
%%%%%%%%%%%%%%%%%%%%%%%% Corollary 5.5 %%%%%%%%%%%%%%%%%%%%%%%%%%%%
\begin{corollary}{\rm (cf. \cite[Corollary 3.3.8]{13})}\label{sec:3.6}
Let $P$ be a non $t$-invertible $t$-prime ideal of a $PVMD$ $R$. Then:
\begin{enumerate}
    \item $P^{-1}=T(P)$ or $T(P)=\Omega(P)$;
    \item If $P\neq (P^2)_t$, then $T(P)=\Omega(P)$;
    \item If $P=(P^2)_t$, then $P^{-1}=T(P)$;
    \item If $P$ is unbranched, then $P^{-1}=T(P)=\Omega(P)$.
\end{enumerate}
\end{corollary}

\begin{proof} \i Follows from Theorem~\ref{sec:3.2}.\\

\ii If $P\neq (P^2)_t$, then there is a prime ideal $Q$ of $R$ such that $\bigcap (P^n)_tR_P=QR_P$. Note that $P\nsubseteq Q$ (otherwise, if
$P=Q$, then $PR_P=QR_P$. But $QR_P\subseteq(P^2)_tR_P=P^2R_P\subsetneqq PR_P$, a contradiction). Hence $T(P)\supseteq R_Q\cap(\bigcap R_{M_{\beta}})\supseteq \Omega(P)$, where $\{M_{\beta}\}$ is the set of all $t$-maximal ideals of $R$ that do not contain $I$. Since $T(P)\subseteq \Omega(P)$, $T(P)=\Omega(P)$.\\

\iii If $P=(P^2)_t$, then $P=(P^n)_t$ for each $n\geq 1$. Hence $(R:P^n)=(R:(P^n)_t)=(R:P)$. So $T(P)=P^{-1}$ by the definition of $T(P)$.\\

\iv Since $P$ is unbranched and $(P^2)_t$ is a $P$-primary (\cite[Proposition 1.3]{19}), $P=(P^2)_t$.
    Hence $T(P)=P^{-1}$ by \iii. It is clear that $\Omega(P)\supseteq T(P)$. By
    \cite[Proposition 1.2]{9}, $P=\bigcup P_{\gamma}$ where
    $\{P_{\gamma}\}$ is the set of primes ideal of $R$ properly
    contained in $P$, and we may assume that they are maximal with this
    property. Then by \cite[Corollary 4.6]{16}, $R_P=\bigcap
    R_{P_{\gamma}}$. Hence by \cite[Theorem 3.2.2]{13}, $\Omega(P)\subseteq
    R_P$. Since $\Omega(P)\subseteq \bigcap R_{M_{\beta}}$, $\Omega(P)\cap R_P\subseteq \bigcap R_{M_{\beta}}\cap R_P$.
    It follows that $\Omega(P)\subseteq P^{-1}=T(P)$.
    Therefore $T(P)=\Omega(P)$.
\end{proof}
\bigskip
%%%%%%%%%%%%%%%%%%%%%%%%%%%%%%%%%%%%%%%%%%%%%%%%%%%%%%%%%%%%%%%%%%%
%%%%%%%%%%%%%%%%%%%% Proposition 5.6 %%%%%%%%%%%%%%%%%%%%%%%%%%%%%%
Our last theorem generalizes \cite[Theorem 3.3.10]{13}.\\
\bigskip
%%%%%%%%%%%%%%%%%%%%%%%%%%%%%%%%%%%%%%%%%%%%%%%%%%%%%%%%%%%%%%%%%%%%%%%%%%%%%%%
%%%%%%%%%%%%%%%%%%%%%%%%%%%%%%%%%%%%%%%%%%%%%%%%%%%%%%%%%%%%%%%%%%%%%%%%%%%%%%%
\begin{thm}\label{sec:3.7}
Let $R$ be a $PVMD$ and $P$ a $t$-prime ideal of $R$. Then:\\
\1 $T(P)\subsetneqq\Omega(P)$ if and only if $T(P)=R_P\cap\Omega(P)$ and $\Omega(P)\nsubseteq R_P$.\\
\2 The following conditions are equivalent:\\
\i $(P\Omega(P))_{t_1}=\Omega(P)$;\\
\ii $\Omega(P)\nsubseteq R_P$;\\
\iii $P=\sqrt{I}$ for some $t$-invertible ideal.
\end{thm}
%%%%%%%%%%%%%%%%%%%%%%%%%%%%%%%%%%%%%%%%%%%%%%%%%%%%%%%%%%%%%%%%%%%%%%%%%%%%%%%%%%%%%%%%%%%%%%%%%%%%%%%%
The proof of this theorem involves the following lemmas. First we notice that in \cite{16}, Gilmer mentioned that $IT(I)=T(I)$ for any invertible ideal $I$ of an arbitrary domain $R$. Our first lemma provides a $t$-analogue result in the class of $PVMDs$. Note that one can replace
the condition ``$PVMD$" on $R$ by assuming that $T(I)$ is a $t$-flat overring of $R$.
\bigskip
%%%%%%%%%%%%%%%%%%%%%%%%%%%%%%%%%%%%%%%%%%%%%%%%%%%%%%%%%%%%%%%%%%%%%%%%%%%%%%%%%%%%%%%%%%%%%%%%%%%%%%%
%%%%%%%%%%%%%%%%%%%%%% Lemma 5.7 %%%%%%%%%%%%%%%%%%%%%%%%%%%%%%%%%
%%%%%%%%%%%%%%%%%%%%%%%%%%%%%%%%%%%%%%%%%%%%%%%%%%%%%%%%%%%%%%%%%%
\begin{lemma}\label{sec:3.8}
Let $I$ be an ideal of a domain $R$.\\
\i If $I$ is $t$-invertible and R is a $PVMD$, then
$(IT(I))_{t_1}=T(I)$ where $t_1$
is the $t$-operation w.r. to $T(I)$.\\
\ii If $I$ and $J$ are two ideals of a domain $R$ such that
$\sqrt{I}=\sqrt{J}$, then $\Omega(I)=\Omega(J)$.
\end{lemma}

\begin{proof}
\i Since $I$ is $t$-invertible, then there is a finitely generated ideal $A$ of
$R$ such that $A\subseteq I_t$ and $A_t=I_t$. Then
$T(I)=T(I_t)=T(A_t)=T(A)=\Omega(I)$ and hence $T(I)$ is a $t$-linked
overring of $R$. Since $I$ is $t$-invertible, then $(II^{-1})_t=R$ and hence
$(I(R:I^n))_t=(R:I^{n-1})$ for each $n\geq 2$. Since $I(R:I^n)\subseteq (I(R:I^n))T(I)$ for
each $n$, then $(I(R:I^n))_t\subseteq(IT(I))_{t_1}$ for each $n$ (
Lemma~\ref{sec:2.6}). Hence $\bigcup(I(R:I^n))_t\subseteq((I(R:I^n))T(I))_{t_1}=(IT(I))_{t_1}$.
So $T(I)=\bigcup(I(R:I^n))_t\subseteq (IT(I))_{t_1}\subseteq T(I)$ and therefore $(IT)_{t_1}=T$, as desired.\\
\ii Straightforward via \cite[Theorem 3.2.2]{13}.
\end{proof}
\bigskip
%%%%%%%%%%%%%%%%%%%%%%%%%%%%%%%%%%%%%%%%%%%%%%%%%%%%%%%%%%%%%%%%%%%
%%%%%%%%%%%%%%%%%%%%%%%%%% Lemma 5.8 %%%%%%%%%%%%%%%%%%%%%%%%%%%%%%
\begin{lemma}{\rm (cf. \cite[Proposition 25.4]{16})}\label{sec:3.9}
Let $R$ be a $PVMD$ and $A_1,....,A_n,B,C\in F_t(R)$. Then:\\
\1 If for each $i$, $A_i$ is $t$-finite, then $\bigcap^{n}_{i=1}
A_i$ is $t$-finite.\\
\2 If $B$ is $t$-finite, then $(C:B)=(CB^{-1})_t$.\\
\3 If $B$ and $C$ are $t$-finite, then $(C:_RB)$ is $t$-finite.
\end{lemma}

\begin{proof}
\1 It suffices to prove it for $n=2$. We have $((A_1\cap A_2)(A_1+A_2))_t=(A_1A_2)_t$ (
\cite[Theorem 5]{17}). Since $A_1$ and $A_2$ are $t$-invertible,
$A_1A_2$ is $t$-invertible and therefore $A_1\cap A_2$ is $t$-invertible and so
$t$-finite.\\
\2 If $x\in (R:B)C$, then $x=\sum_{i=1}^{n}b_ic_i$ where
$b_iB\subseteq R$ and $c_i\in C$. Hence $xB=\sum c_ib_iB\subseteq
RC\subseteq C$. So $(R:B)C\subseteq (C:B)$. Therefore $((R:B)C)_t\subseteq
(C:B)_t=(C:B)$. Conversely, we have $B(C:B)\subseteq C$. Then
$(C:B)=(C:B)_t=((C:B)BB^{-1})_t\subseteq (CB^{-1})_t$.\\
\3  By definition, $(C:_RB)=(C:_RB)_t=((C:B)\cap
R)_t=((CB^{-1})_t\cap R)_t$. Since $C$ and $B$ are $t$-finite,
$(CB^{-1})_t$ is $t$-finite. So by \1, $(C:_RB)$ is $t$-finite.
\end{proof}
\bigskip
%%%%%%%%%%%%%%%%%%%%%%%%%%%%%%%%%%%%%%%%%%%%%%%%%%%%%%%%%%%%%%%%%%%
%%%%%%%%%%%%%%%%%%%%%%% Proposition 5.9 %%%%%%%%%%%%%%%%%%%%%%%%%%%
%%%%%%%%%%%%%%%%%%%%%%%%%%%%%%%%%%%%%%%%%%%%%%%%%%
\noindent{\bf Proof of Theorem 3.6} \1 Assume that $T(P)\subsetneqq\Omega(P)$. Then
$P$ is a non-$t$-invertible $t$-prime  ideal of $R$ (otherwise, if $P$ is
$t$-invertible, then $P$ is $t$-finite, i.e. there is a finitely generated ideal
$A$ of $R$ such that $P=A_t$. Hence $\Omega(P)=\Omega(A_t)=\Omega(A)=T(A)=T(A_t)=T(P)$, a contradiction).
Since $T(P)$ is a subintersection of a $PVMD$ $R$, it is $t$-linked over $R$ (\cite[Theorem 2.22]{25}), and  so $t$-flat over $R$ (\cite[Theorem 2.10]{26}). By Theorem~\ref{sec:3.2}, $P^{-1}=T(P)$. Hence $T(P)=R_P\cap\Omega(P)$ by
\cite[Proposition 1.1]{19} and \cite[Theorem 3.2.2]{13}. Therefore $\Omega(P)\nsubseteq R_P$.\\
The converse is trivial.\\
\2 $(iii)\Rightarrow (i)$ Since $P=\sqrt{I}$,
$\Omega(P)=\Omega(I)$ by Lemma~\ref{sec:3.8}(ii). Since $I$ is $t$-invertible,
by Lemma~\ref{sec:3.8} $(IT(I))_{t_1}=T(I)$. Also since $I$ is
$t$-invertible, there is a finitely generated ideal $A$ of $R$ such that
$A\subseteq I$ and $I_t=A_t$. Hence
$T(I)=T(I_t)=T(A_t)=T(A)=\Omega(A)=\Omega(A_t)=\Omega(I_t)=\Omega(I)$
by \cite[Proposition 3.4]{10}. So
$\Omega(P)=\Omega(I)=(I\Omega(I))_{t_1}\subseteq(P\Omega(I))_{t_1}=(P\Omega(P))_{t_1}\subseteq \Omega(P)$.\\

$(i)\Rightarrow (ii)$ By \cite[Theorem 3.2.2]{13} and \cite[Proposition 4]{5}, $\Omega(P)$ is a $t$-linked overring of $R$. Since $\Omega(P)\nsubseteq R_P$, $(P\Omega(P))_{t_1}=\Omega(P)$ by Lemma~\ref{sec:2.7}(ii).\\

$(ii)\Rightarrow (iii)$ Let $\{Q_{\alpha}\}$ be the set of all $t$-prime ideals
of $R$ that do not contain $P$. Choose $x\in\Omega(P)\setminus R_P$.
Write $x=\frac{a}{b}$ where $a, b\in R$. If $I=(bR:_RaR)$, then
$I\nsubseteq Q_{\alpha}$ for each $\alpha$ and $I\subseteq P$. By
Lemma~\ref{sec:3.9}, $I$ is t-finite and $\sqrt{I}=P$. For this if
$z\not\in\sqrt{I}$, then $z^n\not\in A_v$ for each finitely generated ideal $A$ of $R$ such that
$A\subseteq I$. Hence $z^nab^{-1}\not\in R$ for each $n$. Since
$ab^{-1}\in\Omega(P)$, $z\not\in P$.\quad

\bigskip
%%%%%%%%%%%%%%%%%%%%%%%%%%%%%%%%%%%%%%%%%%%%%%%%%%%%%%%%%%%%%%%%%%%
%%%%%%%%%%%%%%%%%%% Corollary 5.10 %%%%%%%%%%%%%%%%%%%%%%%%%%%%%%%%
\begin{corollary}{\rm (cf. \cite[Corollary 3.3.11]{13})}\label{sec:3.10}
Let $R$ be a $PVMD$ and $P$ a non-$t$-maximal $t$-prime ideal of $R$. Then
$T(P)\subsetneqq\Omega(P)$ if and only if $P=(P^2)_t$ and
$P=\sqrt{I}$ for some $t$-invertible ideal $I$ of $R$
\end{corollary}
%%%%%%%%%%%%%%%%%%%%%%%%%%%%%%%%%%%%%%%%%%%%%%%%%%
\begin{proof}
$\Rightarrow)$ Since $T(P)\subsetneqq\Omega(P)$, $P=(P^2)_t$
(Corollary~\ref{sec:3.6}(ii)) and $\Omega(P)\nsubseteq R_P$ (Theorem~\ref{sec:3.7}). Hence there is a t-invertible ideal of $R$ satisfies
$P=\sqrt{I}$ (Theorem~\ref{sec:3.7}).\\

$\Leftarrow)$ $P=(P^2)_t$ implies that $P^{-1}=T(P)$ by Corollary~\ref{sec:3.6}(iii). Since $P=\sqrt{I}$ for some $t$-invertible ideal $I$ of
$R$, $\Omega(P)\nsubseteq R_P$ by Theorem~\ref{sec:3.7}. By
\cite[Theorem 4.5]{20} $P^{-1}=R_P\cap(\bigcap R_{M_{\beta}})$,
where $\{M_{\beta}\}$ is the set of all $t$-maximal ideals of $R$ that do not
contain$P$. By \cite[Theorem 3.2.2]{13},
$T(P)=P^{-1}=R_P\cap\Omega(P)$. By Theorem~\ref{sec:3.7},
$T(P)\subsetneqq\Omega(P)$
\end{proof}
\bigskip
%%%%%%%%%%%%%%%%%%%%%%%%%%%%%%%%%%%%%%%%%%%%%%%%%
\begin{corollary}{\rm (cf. \cite[Corollary 3.3.12]{13})}\label{sec:3.11}
Let $R$ be a $PVMD$ and $P$ a non-$t$-invertible $t$-prime
ideal of $R$. Then:\\
$(PT(P))_{t_1}\neq T(P)$ and $(P\Omega(P))_{t_2}=\Omega(P)$
if and only if $P^{-1}=T(P)\subsetneqq\Omega(P)$ where $t_1$ (resp. $t_2$) is the $t$-operation w.r. to $T(I)$ (resp.
$\Omega(I)$).
\end{corollary}
%%%%%%%%%%%%%%%%%%%%%%%%%%%%%%%%%%%%%%%%%%%%%%%%%%%%%%
\begin{proof}
If $(PT(P))_{t_1}\neq T(P)$ and $(P\Omega(P))_{t_2}=\Omega(P)$, then
clearly $T(P)\subsetneqq\Omega(P)$. Hence $P^{-1}=T(P)$ by Theorem~\ref{sec:3.2}. Conversely, if $P^{-1}=T(P)\subsetneqq\Omega(P)$, then $(PT(P))_{t_1}\neq T(P)$ by Lemma~\ref{sec:3.4}. Moreover $P=\sqrt{I}$ for some $t$-invertible
ideal $I$ of $R$ by Corollary~\ref{sec:3.10}. Therefore
$(P\Omega(P))_{t_2}=\Omega(P)$ by Theorem~\ref{sec:3.7}.
\end{proof}
%%%%%%%%%%%%%%%%%%%%%%%%%%%%%%%%%%%%%%%%%%%%%%%%%%%%%%%%%%%%%%%%%%%
%%%%%%%%%%%%%%%%%%%%%%%%%%%%%%%%%%%%%%%%%%%%%%%%%%%%%%%%%%%%%%%%%%%
%%%%%%%%%%%%%%%%%%%%    Thebibliography     %%%%%%%%%%%%%%%%%%%%%%%
%%%%%%%%%%%%%%%%%%%%%%%%%%%%%%%%%%%%%%%%%%%%%%%%%%%%%%%%%%%%%%%%%%%

%%%%%%%%%%%%%%%%%%%%%%%%%%%%%%%%%%%%%%%%%%%%%%%%%%%%%%%%%%%%%%%%%%%%%%%%%%%%%%%%%%%%%%%%%%%%%%%%%%%%%%%%%%%%%%%%%%%%
%%%%%%%%%%%%%%%%%%%%%%%%%%%%%%%%%%%%%%%%%%%%%%%%%%%%%%%%%%%%%%%%%%%%%%%%%%%%%%%%%%%%%%%%%%%%%%%%%%%%%%%%%%%%%%%%%%%%%
%%%%%%%%%%%%%%%%%%%%%%%%%%%%%%%%%%%%%%%%%%%%%%%%%%%%%%%%%%%%%%%%%%%%%%%%%%%%%%%%%%%%%%%%%%%%%%%%%%%%%%%%%%%%%%%%%%%%%
\end{document}